\documentclass{amsart}

\usepackage{graphicx}
\usepackage{psfrag}
\usepackage{amsthm}
\usepackage{amscd}
\usepackage{amsfonts}
\usepackage{amstext}

\begin{document}
\title[Complex bounds for multimodal maps]{Complex bounds for multimodal maps: bounded combinatorics}
\author{Daniel Smania }

\address{Instituto de Matem\'atica Pura e Aplicada, 
    Estrada Dona Castorina, 110, 
               Jardim Bot\^anico,
               CEP 22460-320, 
               Rio de Janeiro-RJ,
               Brazil. }
\email{smania@impa.br}

\begin{abstract} We proved the so called complex bounds for multimodal, infinitely renormalizable analytic maps with bounded combinatorics: deep renormalizations have polynomial-like extensions with definite modulus. The complex bounds is the first step to extend the renormalization theory of unimodal maps to multimodal maps. 
\end{abstract}

\maketitle
\newtheorem{lem}{Lemma}[section]
\newtheorem{defi}{Definition}[section]
\newtheorem{rem}{Remark}[section]
\newtheorem{mle}{Main Lemma}[section]
\newtheorem{thm}{Theorem}

\section{Introduction}
\textbf{Renormalization.} A multimodal maps is a smooth map $f\colon I \rightarrow I$ preserving the boundary of $I$ with a finite number of critical points, all of them local maximum or local minimum. Let $f$ be an infinitely renormalizable map in respect to a critical point $p$. This means that there are nested intervals $P^{k}$, $p \in P^{k}$ and  an increasing sequence of natural numbers $N_{k}$ such that 

\begin{itemize}

\item $f^{N_{k}}(P^{k}) \subset P^{k}$.

\item $f^{N_{k}}(\partial P^{k}) \subset \partial P^{k}$.

\item The intervals  $f^{i}(P^{k})$ and  $f^{j}(P^{k})$ have disjoint interiors for  $i \neq j$, $0 \leq i,j < N_{k}$. 

\item The interval $P^{k}$ is maximal with these properties.

\end{itemize}

The map $f^{N_k}$ restricted to $P^{k}$ is the $k$-th renormalization of $f$ with respect to $p$, and it will be called $R^{k}_p(f)$. We say that $f$ has \textbf{bounded combinatorics} with respect to $p$ if $N_{k+1}/N_{k}$ is bounded. The concept of  polynomial-like map was  introduced in \cite{DH}. Sullivan, in \cite{S} showed that the existence of the polynomial-like extensions (for a precise definition of polynomial-like extension, see definition \ref{plike}) of the renormalizations is a useful tool for understanding the dynamics of the infinitely renormalizable maps. We will prove

\begin{thm}[Complex Bounds] \label{cb} Let $f$ be an analytic multimodal map defined in a neighborhood of the interval $I$ such that all critical points have even criticality. Suppose that $f$ is infinitely renormalizable with respect to a critical point $p$, with bounded combinatorics. Then there exists $\epsilon$ such that for sufficiently large $k$, the  $k$-th renormalization with respect to p, $R_{p}^{k}(f)$, has a polynomial-like extension $R_{p}^{k}(f) \colon U \rightarrow V$, with  $mod(V\setminus U) \geq \epsilon$. Furthermore, $diam(V) \leq C |P^{k}|$ and the filled-in Julia set of $R_{p}^{k}(f)$ is contained in a Poincare neighborhood (see \cite{MS}) $D_{\beta}(P^{k})$, where $\beta$ does not depend on $k$.  \end{thm}

Complex bounds, a kind of compactness result, have a lot of applications in the study of the infinitely renormalizable maps: when $f$ is a quadratic-like maps, the proof of the local connectivity of the Julia set involves the complex bounds, and so do the non existence of invariant line fields supported in the Julia set, the convergence of the renormalization operator in the set of infinitely renormalizable maps and the hyperbolicity of this operator in an appropriate space. We expect this theory do be generalized at least for infinitely renormalizable maps with bounded combinatorics.

Complex bounds type results have now a long history. This kind of results was introduced by Sullivan (\cite{S}), in the study of the renormalization operator for unimodal maps with bounded combinatorics. Nowadays, there are a large number of related results. For instance: for infinitely renormalizable unimodal maps  with unbounded  combinatorics, there are independent results by M. Lyubich (\cite{Lyu}) and M. Lyubich \& M. Yampolsky (\cite{LYA97}), by J. Graczyk and G. Swiatek  (\cite{GS96}),  and G. Levin and S. van Strien (\cite{LS}). In the thesis \cite{H95}, there is a proof for infinitely renormalizable bimodal (two critical points) maps with bounded combinatorics in the Epstein class, using the sector lemma (introduced by Sullivan: see \cite{MS} for details), but the proof seems to be incomplete. We will use the methods introduced in \cite{LYA97}. This paper is part of my thesis at IMPA.

\section{Topological Results}

 Assume that all intervals under consideration are closed.

 A multimodal map $f$ is renormalizable with respect to a critical point $p$ if there exist $n > 1$ and an interval $P^{1}$ such that $P^{1}, f(P^1), \dots f^{n-1}(P^1)$ have disjoint interiors and $f^{n}(P^1) \subset P^1$. The smallest $n$ with this property is the period of the renormalization. If there is not an interval containing  properly $P^1$ with the same property, then $f^n(\partial P^1) \subset \partial P^1$. The map $f^{n}$ restricted to $P^1$ is again a multimodal map. This map is called a $renormalization$ of $f$.  Such renormalization could be renormalizable with respect to $p$ again and so on. If this process never finish we say that $f$ is \textbf{infinitely renormalizable} with respect to $p$. Hence we can construct the intervals $P^{k}$ as in the introduction.

 Let $\mathbb{I}_p^{k} = \{ c \colon  f'(c)=0,$ $ c \in \cup_{i < N_{k}} f^{i}(P^{k})   \}$. Since $\mathbb{I}_p^{k+1} \subset \mathbb{I}_p^{k}$, $\mathbb{I}_p^{k}$ is constant for large $k$. Denote $\mathbb{I}_{p} = \mathbb{I}_{p}^{k}$, where $k$ is large, the set of critical points \textbf{involved} in the renormalizations with respect to $p$. 

Assume that there are no wandering intervals for $f$ (for example, if the critical points are non-flat and $f$ is $C^{2}$: see pg. 267 in \cite{MS}). It is well-known that $max \{ |f^{i}(P^{k})|\}_{i < N_{k}}$ goes to zero. In particular for large $k$, there is at most one critical point in each $f^{i}(P^{k})$. We say that an interval $L$ is symmetric with respect to a critical point $c$ if $f$ is monotone in each connect component $J$ of $L \setminus \{c\}$ and $f(J)=f(L)$. Let $J$ be an small interval which contains only one critical point $c$. The symmetrization of $J$ is an interval $L \supset J$ which is a symmetric interval with respect to $c$ and $f(I)=f(L)$. For large $k$,  $P^{k}$ is symmetric with respect to $p$. 

For large $k$, if $ q \in f^{j}(P^{k})\cap \mathbb{I}_p$ and  $r \in f^{j+i \text{ mod }  N_{k}}(P^{k})\cap \mathbb{I}_p$, where $i > 0$ is the smallest number such that $f^{j+i \text{ mod } N_{k}}(P^{k})\cap \mathbb{I}_p \neq \phi$, we say that $r$ is the \textbf{successor} of $q$ at level $k$ and $q$ is the \textbf{predecessor} of $r$ at level $k$.  Denote by $q'_k$ the successor of $q$ at level $k$. 

For $q \in f^{j}(P^{k}) \cap \mathbb{I}_{p}$ denote by the corresponding capital letter $Q_{0}^{k}$ the symmetrization of $f^{j}(P^{k})$. Let $n^{r}_{k}= j - i \text{ mod  } N_{k} $, with $r \in f^{j}(P_{0}^{k})$ and $q \in f^{i}(P_{0}^{k})$, where $q$ is the predecessor of $r$ at level $k$. Then $\{ f^{i}(Q_{0}^{k}) \}_{i < n^{q'_k}_k\text{, } q \in \mathbb{I}_{p}}$ is a family of intervals with disjoint interior. Note that, for large $k$, the boundary of $P^{k}$ contains a periodic point and its symmetric with respect to $p$. Moreover, the orbit of this periodic point does not contain critical points. This implies that $f^{n^{r}_{k}}(Q^{k}_{0}) \subset R_{0}^{k}$, where $r$ is the successor of $q$ at level $k$, and that the  boundary of $Q^{k}_{0}$ contains a point of this periodic orbit.

Note that to prove the  complex bounds (see the precise statement in Theorem \ref{cb} ) for $f$ with respect to $p$, it is sufficient prove the complex bounds for a renormalization with respect to $p$. Replacing $f$ for a deep renormalization with respect to $p$, we can assume, without loss of generality, that $f$ satisfies the \textbf{standard conditions}:

\begin{itemize}
\item $f$ is a composition of unimodal maps with non-flat critical points and  $\{ f(c)\colon f'(c)=0   \}=\{ f(q)   \}_{q \in \mathbb{I}_{p}}$.

\item The interval $Q^{k}_{0}$ is symmetric  with respect to $q$.

\item For large $k$  and for all $q \in \mathbb{I}_{p}$, $Q_{0}^{k}$ contains exactly one critical point.

\item For large $k$, there exists a periodic orbit such that $Q^{k}_{0}$ contains in its boundary a point of this orbit, for all $q \in \mathbb{I}_{p}$.

\item For every $q$, $f^{n^{r}_{k}}(Q^{k}_{0}) \subset R_{0}^{k}$, where $r$ is the successor of $q$ at level $k$. Moreover $r \in f^{n^{r}_{k}}(Q_{0}^{k})$.

\item The intervals in $\{ f^{i}(Q_{0}^{k}) \}_{i < n^{q'_k}_k\text{, } q \in \mathbb{I}_{p}}$  have disjoint interior.
\end{itemize}

\begin{figure}
\centering
\psfrag{d}[][][.8]{$f$}
\psfrag{b}[][][.8]{$f$}
\psfrag{i4}[][][.8]{$P_{0}^{k}$}
\psfrag{i1}[][][.8]{$R_{-(n_{k}^{r}-1)}^{k}$}
\psfrag{i3}[][][.8]{$Q_{-(n_{k}^{q}-1)}^{k}$}
\psfrag{i2}[][][.8]{$Q_{0}^{k}$}
\psfrag{q}[][][.8]{$q$}
\psfrag{p}[][][.8]{$p$}
\psfrag{a}[][][.8]{$f^{n_{k}^{p}-1}$}
\psfrag{c}[][][.8]{$f^{n_{k}^{q}-1}$}
\psfrag{fq}[][][.8]{$f(q)$}
\psfrag{fp}[][][.8]{$f(p)$}
\psfrag{fr}[][][.8]{$f(r)$}
\psfrag{r}[][][.8]{$r$}
\psfrag{e}[][][.8]{$f^{n_{k}^{r}-1}$}
\psfrag{i5}[][][.8]{$R_{0}^{k}$}
\psfrag{i6}[][][.8]{$P_{-(n_{k}^{p}-1)}^{k}$}
\includegraphics[width=0.60\textwidth]{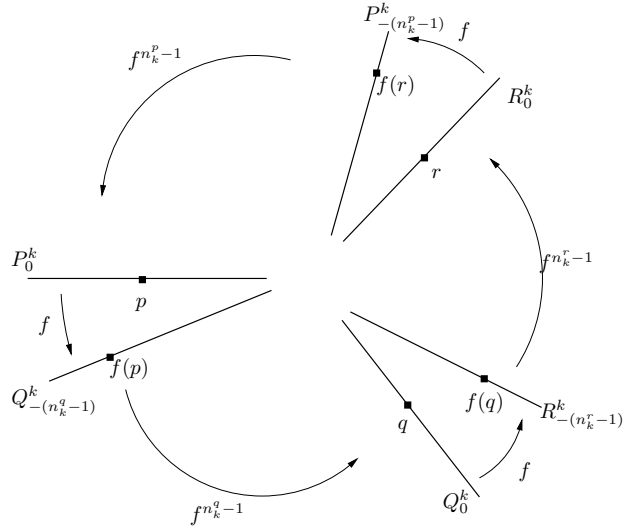}
\caption{Intervals of the  $k$-th renormalization, with 3 critical points involved.}
\end{figure}

We say that a set $W$ is a \textbf{nice set} if it is an union of closed intervals, each interval contains a unique critical point, the critical points of $f$ are contained in the interior of $W$ and $f^{n}(\partial W) \cap int W = \phi$, for all $n>0$.

Let $W$ be a nice set and $D$ be the set of points not in $W$ but whose positive orbit intersects the interior of $W$. The entry map associates to each $x \in D$ the point $f^{n}(x)$ in the interior of $W$ with smallest $n$. It is easy to see that for each connect component $J$ of $D$, $n$ is constant and the image of $J$ by the first entry map is a connect component of $W$. Moreover $J,f(J),\dots,f^{n-1}(J)$ are disjoint connect components of $D$.

\begin{lem} \label{nice} Let $c$ be a critical point. Then there exists $q \in \mathbb{I}_p$ such that $f(c) = f(q)$ and the  connect component $J$ of $f^{-1}(f(Q_{0}^{k}))$ which contains c is a symmetric interval $J$ with respect to $c$ and $f(J)= f(Q_{0}^{k})$. In particular,  the union of connect components of $f^{-1}(f(\cup_{q \in \mathbb{I}_{p}} Q_{0}^{k}))$ which contains a critical point is a nice set. \end{lem}

\begin{proof} If the lemmas does not holds then there exists  critical points  in $J\setminus \{ c \}$. Let  $\tilde{c}$ be the closest one of $c$.  $f(\tilde{c}) \neq f(q)$, otherwise there will be another critical point between $\tilde{c}$ and $c$. But this implies $f(r) \in f(Q_{0}^{k})$, where $r \in \mathbb{I}_{p}\setminus \{ q \}$, which is absurd.
\end{proof}

\begin{rem} \label{rnice} Denote by $\mathcal{N}^k$ the nice set defined in the previous lemma. Let $ \mathcal{A}^{k} = \{  f^{i}(Q_{0}^{k}) \text{ s.t. } 0 < i \leq n^{q'_k}_{k}, \text{ } q \in \mathbb{I}_{p}  \}$. It is easy to see that  $\cup_{J \in \mathcal{A}_k} J \cap \mathcal{N}^k = \cup_{q \in \mathbb{I}_{p}} f^{n_{k}^{q'_k}}(Q_{0}^{k})$. Hence $f^{n_{k}^{r}-i}(Q_{0}^{k}) \subset R_{-i}^{k}$, $i \leq n_{k}^{r} -1$, where $r$ is the successor of $q$ at level $k$ and  $R_{-i}^{k}$ is a domain of the first entry map to $\mathcal{N}^k$. Clearly $f(R_{-(i+1)}^{k})= R_{-i}^{k}$. When we write $R_{-i}^{k}$, assume that $i \leq n_{k}^{r} - 1$.
\end{rem}

We say that an interval $J$ has the property $*_k$ if  there is at most one interval in $\mathcal{A}^{k}$ in the interior of $J$. Note that the closure of $J$ contains at most 3 intervals in $A^{k}$.

For $x \in I$ the intersection number of $x$ with a family of intervals is the number of intervals in the family which contains $x$. The intersection number of a family of intervals is the maximal number of intersection with a point in $I$.

\begin{lem} Let $J$ be an interval such that $f^{n}$ is monotone in $J$ and $\{ f^{i}(J) \}_{i \leq n}$ is a family of intervals with disjoint interior  such that $f^{i}(J)$ intersects an unique interval in the family $\mathcal{A}^{k}$ which is contained in $f^{i}(J)$. Let $L \supset J$ be such that $f^{n}$ is monotone in $L$ and $f^{n}(L)$ satisfies $*_k$. Then the intersection number of the family of intervals $\{ f^{i}(L)\}_{i \leq n}$ is at most three.
\end{lem}

\begin{proof} If the intersection number is at least 4, the interior of $f^{i}(L)$ will contain at least two intervals $H_1$ and $H_2$ in the family $\mathcal{A}^{k}$, for some $i \leq n$. Since $q'_k \in f^{n^{q'_k}_{k}}(Q_{0}^{k})$, for all $q \in \mathbb{I}_{p}$, the intervals  $f^{n-i}(H_1), f^{n-i}(H_1)$ belong to $\mathcal{A}^{k}$, which is absurd, since $f^{n}(L)$ satisfies $*_k$.
 \end{proof}

The previous lemma will allows us to use the real Koebe lemma (Theorem 3.1 and 3.3 in chapter IV in \cite{MS}) without do considerations about the number of intersections of the families of intervals involved.

\begin{lem} \label{pul} Let $U$ be an interval satisfying  $*_k$ such that $Q_{-i}^{k} \subset U$, with $i < n_{k}^{q} - 1$. Then there exists an interval $\tilde{U}$ which contains  $Q_{-(i+1)}^{k}$ such that the map $f$ is monotone in $\tilde{U}$ and $f(\tilde{U})=U$. Furthermore, $\tilde{U}$ satisfies $*_k$. \end{lem}

\begin{proof} Let $\tilde{U}$ be the maximal interval such that $Q_{-(i+1)}^{k} \subset \tilde{U}$,  $f(\tilde{U}) \subset U$ and $f$ is monotone in $\tilde{U}$. We claim that $f(\tilde{U})= U$. Otherwise there exists a critical point $c$ to $f$ in the boundary of  $\tilde{U}$ such that $f(c)$ is in the interior of $U$. Then, by lemma \ref{nice}, the interior of $f(\tilde{U})$ contains $f(R_{0}^{k})$, $r \in \mathbb{I}_{p}$, in the left or right side of $Q_{-i}^{k}$. This proves the first statement. The second statement is obvious, since $f$ restricted to $\tilde{U}$ is monotone and $r'_k \in f^{n_{k}^{r'_k}}(R_{0}^{k})$, for all $r \in \mathbb{I}_{p}$. \end{proof}
 
\begin{lem} \label{pulq} Let $U$ be an interval satisfying $*_k$ which contains $R_{-(n_{k}^{r}-1)}^{k}$, where $q$ is the predecessor of $r$ at level $k$. Then there exists $\tilde{U}$ such that $Q_{0}^{k} \subset \tilde{U}$ and each connect component of $\tilde{U} \setminus Q_{0}^{k}$ is mapped monotonically to the  connect component of $U \setminus R_{-(n_{k}-1)}^{k}$ which contains $f(\partial Q^{k}_{0})$. Furthermore, $\tilde{U}$ satisfies $*_k$. 

\end{lem}

\begin{proof} Similar to previous lemma. \end{proof} 

\begin{defi} \label{pak} Let $U$ be an interval such that  $Q_{0}^{k} \subset U$. The \textbf{pullback of $U$ along the $k$-cycle of renormalization} exists if there are intervals ${\tilde{U}}_{q}$, $U_{r}$, $V_{r}$ for each $r \in \mathbb{I}_{p}$, such that:

\begin{itemize}

\item $U_{q}=U$. For $r \in \mathbb{I}_{p}$, $R_{-(n_{k}^{r}-1)}^{k} \subset V_{r}$ and $f^{n_{k}^{r}-1}(V_{r})=U_{r}$. Furthermore, the map $f^{n_{k}^{r}-1}$ is monotone in $V_{r}$. 

\item  The interval $U_{r}$ is the  symmetric interval with respect to $r$ which contains  $R^{k}_{0}$ and each connect component of $U_{r}\setminus R_{0}^{k}$, is mapped monotonically in a connect component of $V_{c}\setminus C_{-(n_{k}^{c}-1)}^{k}$. Here $c$ is the successor of $r$ at level $k$.

\item The interval ${\tilde{U}}_{q}$ is  the symmetric interval with respect to $q$ which contains  $Q^{k}_{0}$ such that each connect component  of ${\tilde{U}}_{q} \setminus Q_{0}^{k}$ is mapped monotonically in a connect component of $V_{c}\setminus C_{-(n_{k}^{c}-1)}^{k}$. Here $c$ is the successor of $q$ at level $k$.
\end{itemize}
\end{defi}

 Let $M^{k}_{q}$ be the maximal interval which contains $Q_{0}^{k}$ satisfying the property $*_k$ . 

\begin{lem} \label{pula} If $U$ is an interval satisfying $Q_{0}^{k} \subset U \subset M^{k}_{q}$, then the pullback of $U$ along the $k$-cycle of renormalization exists. Furthermore, ${\tilde{U}}_{q}  \subset M^{k}_{q}$, where ${\tilde{U}}_{q}$ is as in the previous definition.
\end{lem}
  
\begin{proof} Immediate consequence of  the previous lemmas. \end{proof}

The $\delta$-neighborhood of an interval $J$ is $\delta$-$J$ $=\{x\colon dist(x,J) \leq \delta |J| \}$. Here $dist(x,J) = \inf \{ |x - y| \colon y \in J   \}$. The mirror image of interval $J$ near a critical point $q$, which is not in $J$, is the other interval $\tilde{J}$ near $q$ such that $f(J)=f(\tilde{J})$.

\begin{rem} \label{rema} In the previous lemma, suppose $\epsilon$-$Q_{0}^{k} \subset U$. Then, by the real Koebe lemma and the quasi symmetry in the critical points, one gets $\delta(\epsilon)-Q_{0}^{k} \subset {\tilde{U}}_{q}$. Furthermore ${\tilde{U}}_{q} \subset Q_{0}^{k-1}$, since the critical points to $f^{N_{k}}$ in ${\tilde{U}}_{q}$ are inside  $Q_{0}^{k}$ and $Q_{0}^{k-1} \setminus Q_{0}^{k}$ contains in both connect components intervals in the form $f^{i}(R_{0}^{k})$, $r \in \mathbb{I}_{p}$, or mirror images of these intervals. Each one of these intervals contains at least one critical point for $f^{N_{k}}$. 
\end{rem}

\begin{lem} Let  $s_{k} = sup\{ |Q^{k}_{-i}|$  s.t. $q \in \mathbb{I}_{p} \}$. Then $s_{k}$ tends to 0. \end{lem}

\begin{proof}
This follows of the non existence of wandering intervals.
\end{proof}

\section{Real Bounds} \label{sect}
An interval $L$ is  $\delta$-small with respect to level $j$ if $L \subset Q_{0}^{j-1}$ and $|L|/|Q_{0}^{j}| \leq \delta$, for some $q \in \mathbb{I}_{p}$. The interval $J_1$ is $\delta$-deeply contained in $J_2$, $ J_1\subset_{\delta} J_2$, if  $\delta$-$J_1 \subset J_2$. The intervals $J_1$ and $J_2$ are $C$-commensurable if $1/C \leq |J_1|/|J_2| \leq C$. Let $P(f)$ be the poscritical set  $\{f^{i}(c)\colon f'(c)= 0  \}$. Denote by $Fol_{\ell}(x) = |x|^{\ell}$ the folding function with criticality $\ell$. Let $A_{J}\colon J \rightarrow I$ be an affine map which maps the interval $J$ in $I=[-1,1]$. When we write $J=[a,b]$, $a \neq b$, we do not assume $ a < b$.

We are going to prove real estimates. In \cite{H95} and \cite{HU98}, J. Hu obtains real estimatives for infinitely renormalizable multimodal maps (with  weaker smooth assumptions), similar to lemmas \ref{bg} and \ref{realb}. However, to use Lyubich-Yampolsky approach we need other estimates (see lemma \ref{ts}). For completeness we present the full argument. In this section $f$ is a $C^3$ multimodal map satisfying the standard conditions.

\begin{defi} One of the points in the boundary of $Q_{-i}^k$ is periodic and the other is eventually periodic. Let $A_{Q_{-i}^k}\colon Q_{-i}^k \rightarrow I$ be the affine map which maps the periodic point in  -1. For $q \in \mathbb{I}_p$, define $\phi_q^{k}\colon I \rightarrow I$ and $\psi_{q}^k \colon I \rightarrow I$ by
\begin{equation}
\phi_q^{k} = A_{Q_{0}^{k}}\circ f^{n_{k}^{q}-1} \circ A^{-1}_{Q_{-(n_{k}^{q} -1) }^{k}}
\end{equation}

\begin{equation}
\psi_{q}^k = A_{R_{-(n_{k}^{r} -1) }^{k}} \circ f \circ A^{-1}_{Q_{0}^{k}}.
\end{equation}

 Here the critical point $r$ is the successor of $q$ at level $k$.

 The maps $\phi_{q}^{k}$ are the \textbf{monotone parts} of the $k$-th renormalization. The \textbf{partial} monotone parts are the diffeomorphisms  $ A_{Q_{0}^{k}} \circ f \circ A_{Q_{-i }^{k}}^{-1}$. The maps $\psi_{q}^k$ are the \textbf{folding parts} of the renormalization.    

\end{defi}

Let us to begin with an obvious lemma:

\begin{lem}\label{unif} Consider $g = h_1 \circ Fol_\ell \circ h_2$, where $h_1$, $h_2$ are $C^1$ diffeomorphisms in a neighborhood of zero.  There exists $C$ such that if $J$, $M$ are sufficiently small intervals such that $J$ contains the critical point $c$ of $g$ and  $g(J)\subset M$ then $|D(A_M \circ g \circ A^{-1}_J)(x)| \leq C |x-\tilde{c}|^{\ell-1}$. Here $\tilde{c}$ is the critical point of $A_M \circ g \circ A^{-1}_J$.
\end{lem}
\begin{proof} Use that $h_1$ and $h_2$ are bi-Lipchitz in a neighborhood of zero.
\end{proof}

The following lemma will be used full-time, sometimes without explicit mention:

\begin{lem}[Theorem 8.1 in \cite{NS}]\label{atra} For all $l > 1$, $C_1, C_2 > 0$, and $n_0 \in \mathbb{N}$, there exists $\delta$ such that if 

\begin{itemize}
\item The multimodal map $g\colon J \rightarrow J$ satisfies $max_J |Dg| \leq C_1$.
\item  For a critical point $c$ of $g$,  $|Dg(x)| \leq C_2 |x-c|^{\ell-1}$, for all $x \in J$.
\item  There is $x_0 \in J$ with $|x_0 - c| \leq \delta$ and $| f^{n}(x_0) - c| \leq \delta$, $n \leq n_0$.
\end{itemize}
Then there is an hyperbolic attractor which contains $c$ in its immediate basin.
\end{lem}

Note that critical points in $\mathbb{I}_p$ can not be attracted by a periodic orbit, because $\mathbb{I}_p=\mathbb{I}^{k}_p$ for large $k$.

\begin{lem}[Real Bounds: see \cite{H95}, \cite{HU98}] \label{realb} For large $k$, the following holds:

\begin{itemize}

\item There exists $C > 0$ such that  $1/C \leq |D \phi_q^{k}| \leq C$, $q \in \mathbb{I}_p$. The same control of derivatives holds for the partial monotone parts. 

\item There exists $K > 0$ such  that $|D \psi_{q}^k(x)| \leq K |x-c|^{\ell-1}$, where $c$ is the critical point of $\psi_{q}^k$ and $\ell$ is the criticality of $f$ in $q$.

\item There are $\lambda_{1}$ and $\lambda_{2}$ such that $0 < \lambda_{1} < |Q_{0}^{k+1}|/|Q_{0}^{k}| < \lambda_{2} < 1$.
\end{itemize}
\end{lem}

\begin{proof}
We will use the smallest interval trick. Let $Z = f^{i}(C^{k}_{0})$ be the smallest interval in $\{  f^{i}(C_{0}^{k}) \colon c \in \mathbb{I}_{p} \text{, } 0 \leq i < n_{k}^{c'_k} \} $. 
Observe that  $\delta$-$Z \subset I$, where $\delta$ depends on $max f'$. Taking $\delta \leq 1/2$, the interval $\delta$-$Z$ satisfies $*_k$. Consider the maximal interval $X$ which contains $f(C^{k}_{0})$ such that $f^{i-1}$ is monotone in $X$ and $f^{i-1}(X) \subset$ $\delta$-$Z$. We claim that $f^{i-1}(X)=$ $\delta$-$Z$. Otherwise, there exists $a \in \partial X$ such that $f^{j}(a)$ is a critical point for $f$, $ 0 < j < i - 2$. By lemma \ref{nice} and remark \ref{rema},  $f^{j+1}(X) \setminus f^{j + 2}(C^{k}_{0})$ contains an interval in $\mathcal{A}^{k}$. This is a contradiction.
By the real Koebe lemma and the quasi symmetry at the critical points, there exists $\tilde{\delta}$ such that $\tilde{\delta}$-$C^{k}_{0}$ satisfies $*_k$. Use lemmas  \ref{pul} and \ref{pulq} many times and the real Koebe lemma to conclude that there is $\epsilon$ such that $\epsilon$-$Q^{k}_{0}$ satisfies $*_k$. In particular, $M^{k}_{q}$ is a $\epsilon$-neighborhood of $Q^{k}_{0}$. By lemma \ref{pula}, real Koebe lemma and the quasi symmetry at the critical points, one gets control of distortion in the monotone parts of  the k-cycle and that the pullback of $M^{k}$ along the k-cycle is a $\delta(\epsilon)$-neighborhood of $Q^{k}_{0}$. Moreover, this neighborhood is contained in $Q^{k-1}_{0}$, by remark \ref{rema}. Hence  $|Q_{0}^{k}|/|Q_{0}^{k-1}| < \lambda_{2} < 1$. Since $|C^{k}_{0}|$ goes to zero, by lemma \ref{unif}, one gets the second statement.  In particular,  if $Q_{0}^{k}$ is small with respect to $Q_{0}^{k-1}$ then $Q_{0}^{k}$ contains an attractor which attracts $q$, because the combinatorics is bounded and by the previous lemma. But this is absurd, since $\mathbb{I}_p^{k}$ is constant for large $k$. 
\end{proof}

\begin{lem} \label{Sn} The map $f^{n_k^{q}}$ has negative Schwartzian derivative in $Q_{-(n_{k}^{q}-1)}^{k}$, $q \in \mathbb{I}_{p}$ and large $k$. 
\end{lem}
\begin{proof} Let $C_2 = max Sf$. Near of $q$, we have $Sf(x) \leq  -C_1/(x-q)^{2}$. By the real bounds
\begin{equation}
 Sf^{n_{k}^{q}}(x) = Sf(f^{n_{k}^{q} -1}(x))|Df^{n_{k}^{q} -1 }(x)|^{2} + \sum_{t < n_{k}^{q} - 1} Sf(f^{t}(x))|Df^{t}(x)|^{2}
\end{equation}

\begin{equation}
\leq |Df^{n_{k}^{q} -1 }(x)|^{2}( -\frac{C_1}{|Q_{0}^{k}|^{2}} + \sum_{t < n_{k}^{q} - 1} Sf(f^{t}(x))|Df^{n_{k}^{q} -1 - t}(f^{t}(x))|^{-2})
\end{equation}
\begin{equation}
\leq \frac{|Df^{n_{k}^{q} -1 }(x)|^{2}}{|Q_{0}^{k}|^{2}}( - C_1 + C_3 C_2 \sum_{0 < t \leq n_{k}^{q} - 1} |Q_{-t}^{k}|^{2} )
\end{equation}
\begin{equation}
\leq \frac{|Df^{n_{k}^{q} -1 }(x)|^{2}}{|Q_{0}^{k}|^{2}}( - C_1 + C_3 C_2 max\{ |Q_{-t}^{k}|\}_{0 < t \leq n_{k}^{q} - 1})
\end{equation}

We learn this argument in \cite{Ko97}.
\end{proof}

 \begin{lem} \label{icom} For every  $C_{1} \in \mathbb{N}$, there exists  $C_{2}$ such that , for all  $\ell$ with  $Q_{-\ell}^{k} \subset  R_{0}^{j}$, $j<k$ and  $\ell < C_{1}N_{j}$, we have:

\begin{equation}
 \frac{|R_{0}^{j}|}{|Q_{-\ell}^{k}|}\leq C_{2} \frac{|Q_{0}^{j}|}{|Q_{0}^{k}|}
\end{equation}

\textbf{ $C_{2}$ does not depend on $j$,$\ell$ or $k$}. \end{lem}

\begin{proof} Because $\ell < C_{1}N_{j}$, $Q_{-\ell}^{k}$ will return to $Q_{0}^{k}$ after a bounded number of iterations of monotone and folding parts of the $j$-th renormalization. In particular, by real bounds, the derivative of map
$$A_{Q_{0}^{j}} \circ f^{\ell} \circ A_{R_{0}^{j}}^{-1} $$
is bounded away from infinity. This proves the lemma. 
\end{proof}

\begin{lem}\label{cut} The point $f^{i}(q)$ cuts  $R_{-(n_{k}^{r} -i) }^{k}$ in commensurable parts, where $q$ is the predecessor of $r$ at level $k$.\end{lem}

\begin{proof}  Assume by contradiction that  $f^{i}(q)$ is very close to a point $a \in \partial R_{-(n_{k}^{r} -i) }^{k}$ with  $f^{n_{k}^{r} -i}(a) \in \partial R_{0}^{k}$. By the real bounds for the partial monotone parts of the $k$-th renormalization, $f^{n_{k}^{r}}(q)$ is very close to $\partial R_{0}^{k}$ with respect to level $k$. But $f^{n_{k}^{r}}(q)$ returns to $R_{0}^{k+1}$ after a bounded number of iterates of $f^{N_{k}}$. Since $R_{0}^{k+1}$ is deeply inside $R_{0}^{k}$ and by the  control of derivatives of $f^{N_{k}}$ in $R_{0}^{k}$, this is a contradiction. \end{proof}

\begin{lem} \label{cuts} Let  $J$ be a interval in the family $\{ f^{i}(Q_{0}^{k}) \}_{0 < i < n^{q'_k}_{k}}$. There exists a critical point $c$ for $f^{N_{k}}$ which cuts $J$ in  commensurable parts. \end{lem}

\begin{proof}
Let $f^{j}(q) \in \partial f^{j}(Q_{0}^{k})$. Let $b$ be the closest point of $f^{j}(q)$ in $f^{j}(Q_{0}^{k})$ such that $f^{N_k - j}(b)=q$. Note that if $b$ is close to $f^{j}(q)$, using lemma \ref{realb} and that $f^{N_k}(b)=f^{j}(q)$,  there will be an attractor in the interval $[f^{j}(q),b]$ which attracts a critical point. This is an absurd. The point $b$ cannot be close to the periodic point in $\partial f^{j}(Q_{0}^{k})$, otherwise, by the real bounds, $f^{N_k - j}(b)=q$ would be close to the boundary of $Q_{0}^{k}$, which is impossible. 
\end{proof}

\begin{lem}\label{com} There exists $\delta$ such that if $i \leq n^{q'_k}_k$ and $f^{i}(Q_{0}^{k}) \subset R_{0}^{k-1}$, for some $q,r \in \mathbb{I}_{p}$, then $|f^{i}(Q_{0}^{k})|/|R_{0}^{k-1}| \geq \delta$.
\end{lem}

\begin{proof} Note that  $f^{N_{k}-i}(f^{i}(Q_{0}^{k}))$ is commensurable with $Q_{0}^{k-1}$. By the real bounds, the map
$$ A_{Q_{0}^{k-1}}\circ f^{N_{k}-i} \circ A_{R_{0}^{k-1}}^{-1}$$
has bounded derivative. This proves the lemma.
\end{proof}

\begin{lem} \label{ts} Let $L^{k}_{q}$ be the maximal interval $L$ which contains $Q_{0}^{k}$ and  $f^{N_{k}}$ has no critical points in $L \setminus Q_{0}^{k}$. Then there exists $\delta$, which does not depend on $k$,  and  intervals $S_{q}^{k}$ , $T_{q}^{k} $ such that, for large $k$
\begin{enumerate}

\item These  intervals are nested: $Q_{0}^{k} \subset_{\delta} S_{q}^{k} \subset L^{k}_{q} \subset_{\delta} T_{q}^{k} \subset Q_{0}^{k-1}$ .

\item The interval  $S_{q}^{k}$ is the pullback of  $T_{q}^{k}$ along the k-cycle and  this  pullback can be extended to the $\delta$-neighborhood of $T_{q}^{k}$ which is also contained in $Q_{0}^{k-1}$.
\end{enumerate}
\end{lem}

\begin{proof} Let $M^{k}_{q}$ be as in lemma \ref{pula}. In the proof of lemma \ref{realb} we proved that $M^{k}_{q}$ contains a $\delta$-neighborhood of $Q_{0}^{k}$.  Let $T_{p}^{k} = M^{k}_{q} \cap Q_{0}^{k-1}$. Let $S_{p}^{k}$ be the pullback of $T_{p}^{k}$ along the k-cycle. By remark \ref{rema} and lemma \ref{pula},  $Q_{0}^{k} \subset_{\delta} S_{p}^{k} \subset T_{p}^{k}$. Let $J$ be a connect component of $T_{p}^{k}\setminus Q_{0}^{k}$. There are two cases: if $J$ is a connect component of $Q_{0}^{k-1} \setminus Q_{0}^{k}$ then $J$ clearly contains an interval in the family $\mathcal{A}^k$ or a mirror image of one. Otherwise  $J$ is a connect component of $M^k_q \setminus Q_{0}^{k}$ and it contains an interval in the family $\mathcal{A}^k$, by the maximality of $M^k_q$.  Note that this interval contained in $J$ is commensurable with $Q_{0}^{k}$(lemma \ref{com}) and it contains a critical point which cuts it in commensurable parts(lemma \ref{cuts}).  Hence $L^{k}_{q} \subset_{\delta} T_{q}^{k}$, since $L^{k}_{q} \setminus Q_{0}^{k}$ does not contains critical points for $f^{N_k}$. Replacing $T_{q}^{k}$, if necessary, by a smaller interval between $T_{q}^{k}$ and $L^{k}_{q}$, we obtain the second statement of ii.
\end{proof}

\begin{lem} Let $u$ be the periodic point in $\partial Q_{0}^{k}$. There exist $\delta$ such that if $J$ is the maximal interval of monotonicity of $f^{N_{k}}$ which contains $u$, then $J \subset_{\delta} f^{N_{k}}(J)$. 
\end{lem}

\begin{proof} Note that $J\setminus Q_{0}^{k}$ is the connect component of $L^k_q\setminus Q_{0}^{k}$ which contains $u$. Hence $J$ is commensurable with $Q_{0}^{k}$ and   $J\setminus Q_{0}^{k} \subset_{\delta} f^{N_{k}}(J) \setminus Q_{0}^{k}$, for some $\delta$, by lemma \ref{ts}.i. Let $c$ be the boundary point of $J$ inside $Q_{0}^{k}$. Surely $[u,c] \subset [u.f^{N_k}(c)]$. It is sufficient to prove that the interval $[c,f^{N_k}(c)]$ is commensurable with $Q_0^k$. Indeed, if $[c,f^{N_k}(c)]$ is very small, by the real bounds  and lemma \ref{atra} there will be an attractor which contains a critical point in its immediate basin. This is absurd, since $\mathbb{I}^{k}_p$ is constant for large $k$.
\end{proof}

\begin{lem}[Bounded geometry: see \cite{H95},\cite{HU98}] \label{bg} The geometry of the intervals $Q_{i}^{k}$ is bounded in the following sense:

\begin{itemize}
\item If $Q_{-j}^{k} \subset C_{-l}^{k-1}$, then the distance between $Q_{-j}^{k}$  and $\partial C_{-l}^{k-1}$ is commensurable with $|Q_{-j}^{k}|$.

\item Suppose $R_{-i}^{k},Q_{-j}^{k} \subset C_{-l}^{k-1}$, with $\partial R_{-i}^{k} \cap \partial  Q_{-j}^{k}= \phi$. Then the space between these intervals is commensurable with their length. 

\item In particular, $\sum_{q \in \mathbb{I}_p} \sum_{i < n_q^{k}} |Q_{-i}^{k}|$ goes exponentially fast to zero.
\end{itemize}
\end{lem}

\begin{proof} Observe that $Q_{-j}^{k}$ is not small with respect to $C_{-l}^{k-1}$. Otherwise, by the real bounds and bounded combinatorics, there will be an attractor in $Q_{-j}^{k}$ which attracts the  critical point $q$. Suppose that a point in the boundary of $Q_{-j}^{k}$ is close to the boundary of $C_{-l}^{k-1}$. By the control on the derivatives in the monotone and folding parts of the $(k-1)$-th renormalization, such point will return to $Q_{0}^{k}$ after a large number of iterations of the monotone and folding parts of the $(k-1)$-th renormalization, which is absurd, because the combinatorics is bounded. To  prove the second statement, let $a_{1} \in \partial R_{-i}^{k}$ and $a_{2} \in \partial Q_{-j}^{k}$ be such that $dist(R_{-i}^{k}, Q_{-j}^{k})=dist(a_1,a_2)$. Let $t_1, t_2$ be the minimal times such that $f^{t_1}(R_{-i}^{k}) \subset R_{0}^{k}$ and $f^{t_2}(Q_{-i}^{k}) \subset R_{0}^{k}$. By the  lemma \ref{ts}, the real Koebe lemma and the quasi symmetry at the critical points, for every $D_{-i}^{k}$, $d \in \mathbb{I}_{p}$, there is in each side of $D_{-i}^{k}$  an interval $J$ commensurable with $D_{-i}^{k}$ such that $f^{a}$ is monotone in $J$. Here, $a$ is the time of first entry of $D_{-i}^{k}$ to $R_{0}^{k}$. Hence if the space between $R_{-i}^{k}$ and $Q_{-j}^{k}$ is small and $t=max\{t_1,t_2  \}$, then $f^{t}(a_1)$, $f^{t}(a_2)$ are very close distinct points in $R_{0}^{k-1} \setminus R_{0}^{k}$. Moreover $f^{t+1}(a_1)$, $f^{t+1}(a_2)$ are distinct points in the orbit of the periodic point contained in $\partial R_{0}^{k}$. By the control of derivatives of $(k-1)$-renormalization and bounded combinatorics, the interval $[f^{t+1}{a_1},f^{t+1}{a_2}]$ is mapped monotonically into itself by the map $f^{N_k}$, which is absurd, since the negative Schwartzian derivative of $f^{N_{k-1}}$ in $f(R_{0}^{k-1})$ implies that an attractor in $[f^{t+1}{a_1},f^{t+1}{a_2}]$ contains a critical point in its immediate basin.
\end{proof}

\begin{rem} We can assume $f \in C^{2}$ with non flat critical points, or even weaker smoothness assumptions,  and substitute the lemma \ref{Sn} by a theorem in  pg. 268 of \cite{MS}, which claims that a no repelling periodic orbit with large period must have a critical point in its immediate basin.
\end{rem}

\begin{rem} \label{zeror} Let $U_q=T_{q}^{k}, U_r, V_r, \tilde{U}_q=S_{q}^{k}$ be as in definition \ref{pak}. For each $r \in \mathbb{I}_p$, there is a neighborhood $\tilde{V}_r$ of $V_r$ which is mapped monotonically by $f^{n^r_k -1}$ to a $\delta$-neighborhood of $U_r$, for some $\delta$, by lemma \ref{ts}.ii, the real Koebe lemma and the quasi symmetry in the critical points.  Reducing $\delta$ a little bit we obtain distortion control in $\tilde{V}_r$, by real Koebe lemma. Hence, since $\sum_{i < n^r_k} |R^k_{-i}|$ goes exponentially fast to zero, so do $\sum_{i < n^r_k} |f^{i}(\tilde{V}_r)|$. \end{rem}

\begin{lem} There is $\delta$ such that, for all $c \in \mathbb{I}_{p}$, $P(f) \cap \text{$\delta$-$Q_{0}^{k}$} = P(f)\cap Q_{0}^{k}$.
\end{lem}

\begin{proof}
 If the statement does not holds then there exists $f^{i}(c) \in I\setminus Q_{0}^{k}$, for some $c \in \mathbb{I}_{p}$, very close to $Q_{0}^{k}$ with respect to level $k$, hence $f^{i}(c) \in Q_{0}^{k-1}$. After a bounded number of iterations of $f^{N_{k-1}}$, $f^{i}(c)$ is mapped inside $Q_{0}^{k+1}$. But this is impossible, because $f^{i}(c)$ is very close to a point in $\partial Q_{0}^{k}$, which never returns to $Q_{0}^{k}$.
\end{proof}

\section{Monotone Pullbacks} 
In the rest of this paper, $f$ is a multimodal map satisfying the standard conditions, analytic in a neighborhood of $I$, whose critical points have even criticality. In particular, for each $q \in \mathbb{I}_p$, there are univalent maps $\tilde{h}_q$ and $h_q$ defined in neighborhoods of $f(q)$ and $q$, respectively, such that these maps are real in real points and $f = \tilde{h}^{-1}_q  \circ Fol_\ell \circ h_q$ near to $q$. Here $Fol_{\ell}(x)=x^\ell$, where $\ell$ is even and may vary with $q$. Assume that if $J$ is a small interval such that $q \in J$, then $h_q(J) \subset \mathbb{R}^{+}$. If $g$ is a inverse branch of  $Fol_\ell$, we say that  $h_q \circ g \circ  \tilde{h}^{-1}_q$ is an \textbf{inverse branch} of $f$ near to $f(q)$. Consider the inverse branches $Fol^{+}_{\ell}\colon \mathbb{C}_{\mathbb{R^{+}}}\rightarrow \mathbb{C}$ and $Fol^{-}_{\ell}\colon \mathbb{C}_{\mathbb{R^{+}}}\rightarrow \mathbb{C}$ such that $Fol^{+}_{\ell}(1)=1$ and $Fol^{-}_{\ell}(1)=-1$. Define the \textbf{real inverse branches} of $f$ near to $f(q)$ by 
\begin{eqnarray}
f^{+}_q =   h_q \circ Fol_{\ell}^{+} \circ  \tilde{h}^{-1}_q \\
f^{-}_q =    h_q \circ Fol_{\ell}^{-} \circ  \tilde{h}^{-1}_q
\end{eqnarray}

Fix a level $k_0$ of renormalization deep enough such that for each $r \in \mathbb{I}_p$, $R_{0}^{-(n^{r}_{k_0}-1)}$ is deeply contained in the closure of the domain of $f^{+}_q$ and $f^{-}_q$, for all $r \in \mathbb{I}_p$, where $q$ is the predecessor of $r$ at level $k_0$. Let $\delta_1$, $\delta_2$ be such that

\begin{itemize}
\item  The map $f$ restricted to $B_{\delta_1}(Q^{k_0}_{-i})$ is univalent and  $B_{\delta_2}(Q^{k_0}_{-(i-1)}) \subset f(B_{\delta_1}(Q^{k_0}_{-i}))$, for all $i$ such that  $0 < i \leq n^{q}_{k_0} -1 $ and  $q \in \mathbb{I}_p$
\end{itemize}

We say that the interval $J$ is in the level $k_0$ if $J \subset Q^{k_0}_{-i}$, for some $q \in \mathbb{I}_p$, $0 \leq i \leq n^{q}_{k_0} -1$.

\begin{defi} Let $J$ be an interval in the level $k_0$ such that $f^n$ is monotone in $J$, for some $n$. Let $z \in \mathbb{C}$. We say that  the \textbf{ backward orbit $z_{-i}$ of $z$ along the orbit of $J$ by $f^{n}$} is well defined if for  $i \leq n$, $z_{-i}$ satisfies the following properties:

\begin{itemize}
\item $z_0 = z$. If $f^{n-i}(J) \subset Q^{k_0}_{-l}$, $l <  n^{q}_{k_0} - 1$, then $z_{-i} \in B_{\delta_2}(Q^{k_0}_{-l})$ and $z_{-(i+1)}$ is the unique point in $B_{\delta_1}(Q^{k_0}_{-(l+1)})$ such that $f(z_{-(i+1)})=z_{-i}$.

\item If $f^{n-i}(J) \subset Q^{k_0}_{-(n^{q}_{k_0} - 1)}$, then $z_{-i}$ is contained in the domain of $f^{+}_q$ and $f^{-}_q$. Furthermore 

\begin{equation}
z_{-(i+1)}= \begin{cases}
f^{+}_q(z_{-i})&     \text{if $ f^{+}_q(f^{n-i}(J)) = f^{n-(i+1)}(J)$},   \\
 f^{-}_q(z_{-i})&    \text{if  $f^{-}_q(f^{n-i}(J)) = f^{n-(i+1)}(J)$.}
\end{cases}
\end{equation}
\end{itemize}
We say that $z_{-n}$ has the \textbf{itinerary} of $J$ by $f^{n}$. Let $B \subset \mathbb{C}$ such the backward orbit $z_{-i}$ is well defined for each $z \in B$. Then the set $\{z_{-n}\colon z \in B\}$ is the \textbf{complex pullback}   of $B$ until $J$ by $f^n$. 
\end{defi}

For an interval $J$, let $B_\delta(J)$ be the set of complex numbers whose distance to $J$ is at most $\delta$. Let $\mathbb{C}_{J} = \mathbb{C}\setminus(\mathbb{R}\setminus J) $. We will denote by $D_{\theta}(J)$ the  Poincare's neighborhood of  $J$ with internal angle $\theta$ (see pg. 485 in \cite{MS}, or \cite{S}). Denote by $D(J)$ the Poincare's neighborhood with internal angle $\pi/2$. Given $z \in \mathbb{C}$ and an interval $J=[a_0,a_1]$, let $\widehat{(z,a_i,J)}$ be the angle between the segment $[a_i,z]$ and the real ray beginning at $a_i$ which does not contains $J$. Let $\widehat{(z,J)}=min\{ \widehat{(z,a_1,J)},\widehat{(z,a_2,J)}\}$. The $\epsilon$-sector with vertex $a_i$ is the set of complex numbers $z$ such that $\widehat{(z,a_i,J)} \leq \epsilon$.

Conditions which allow us to make pullbacks  of points along the orbit of an interval will be the main interest in this section. In Epstein maps, the classical tool is the Schwartz lemma:  Poincare neighborhoods  are mapped by Epstein maps to inside Poincare neighborhoods with same internal angle (see \cite{S}). The following results will  allow  us to use the same principle to analytic maps (compare with  \cite{EF}). 

\begin{lem}[Almost preserved Poincare neighborhood] \label{pulb} There exists $K$ with the following property. There exists $l_0(\theta_0)$ such that if $J$ is an interval in the level $k_0$ satisfying

\begin{itemize}
\item $f^{n}$ is monotone in $J$.
\item $\sum_{i \leq n} |f^{i}(J)| \leq l_0$.
\end{itemize}

Then, provided that $\theta \leq \theta_0$,  if $z \in D_{\theta}(f^{n}(J))$,  the backward orbit of $z$ along the orbit of $J$  is well defined and $z_{-n} \in D_{\tilde{\theta}}(J)$, with $\tilde{\theta} = \theta exp(K\sum_{i \leq n} |f^{i}(J)|)$.
\end{lem}

\begin{proof} Since the maps $f^{+,-}_{q}$ are  compositions of univalents maps which map the real line to the real line and the maps $Fol_{\ell}^{+,-}$ belong to the Epstein class, we can use lemma 5.2 (see also 5.3)   in pg. 487 in \cite{MS}.
\end{proof}

\begin{lem}[Universal distortion control] \label{cdi} \label{nodist} For all $\delta_0$ and  $\theta_0$ there exist $C$ and $\epsilon_0$ with the following properties. Let $h \colon  D_{\theta}($$\delta$-$M)\rightarrow \mathbb{C}_{M'}$ be an univalent map, $h(M)=M'$, with $\theta \leq \theta_0$ and $\delta \geq \delta_0$. Then 
\begin{enumerate}
\item  If $z \in D_{\theta}(M)$ and the interval $J$ is contained in $M$ then
\begin{equation}
\frac{dist(h(z),h(J))}{|h(J)|} \leq C\frac{dist(z,J)}{|J|}
\end{equation}

\item Provided $\epsilon < \epsilon_0$, $\widehat{(z,J)}\leq \epsilon$ implies $\widehat{(h(z),h(J))}\leq \alpha=O(\epsilon)$.

\end{enumerate}
\end{lem}

\begin{proof}
Use the control on the first and second derivatives given by complex Koebe lemma.
\end{proof}

The following four lemmas are a straightforward modification of  statements contained in \cite{LYA97}.

\begin{lem} \label{cpulq} Let $C \geq 0$, $\theta_0 > 0$. There exists $l_0(\theta_0)$ with the following property. Let $J$ and $J'$ be intervals such that
\begin{itemize}
\item  $|J| \leq l_0$.
\item  $J$ contains  $f(q)$ and  $1/C \leq |J_1|/|J_2| \leq C$, where $J_1$,$J_2$ are the connect components of $J\setminus \{f(q)\}$.
\item $J'$ is the maximal symmetric interval with respect to $q$ such that $f(J')$ is contained in $J$.
\end{itemize}
If  $\theta < \theta_0$ and  $g$ is an inverse branch of $f$ near to $f(q)$ then $g(D_{\theta}(J)) \subset D_{\tilde{\theta}}(J')$. The angle $\tilde{\theta}$ depends on $\theta$ and $C$. Note that the inverse branch $g$ \textbf{does not need} to be real.
\end{lem}

\begin{lem}\label{naomuda} If $\epsilon$ is small enough, the following holds. Let $J=[a,b]$ be a small interval near to $q$ and  $z \in \mathbb{C}$ near to $q$. Then $\widehat{(f(z),f(a),f(J))}\leq \epsilon$ implies $\widehat{(z,b,J)} > \epsilon$.
\end{lem}

\begin{lem} \label{mi} There exists $l_0$, $C > 0$ and $\epsilon_0$ such that if 
\begin{itemize}
\item The interval $J$ is contained in $L$, which is a maximal interval of monotonicity of $f^{n}$. $L$ is at level $k_0$. 

\item $\sum_{i \leq n}|f^{i}(L)| \leq l_0$.

\item $z$ is a point with the same itinerary of $J$ by $f^{n}$ and $\widehat{(f^i(z),f^i(J))}\leq \epsilon \leq \epsilon_0$, for all $i < n$.
\end{itemize}

Then $z \in D_{\alpha}(L)$, where $\alpha = (\pi/2 + O(\epsilon))exp(C\sum_{i \leq n}|f^{i}(L)|)$.
\end{lem}

\begin{lem} \label{liu} For any $\epsilon$, $\delta$ and  $l_0$, with $l_0$ small enough, there exists $K(\epsilon,\delta,l_0)$ with the following property. If the interval $J$ and $z$ satisfy 

\begin{itemize}
\item There exists an interval $L$ in the level $k_0$ which contains $J$, such that $f^{n}$ is monotone in $L$. Moreover $\delta$-$f^{n}(J) \subset f^{n}(L)$.

\item $\sum_{i \leq n}|f^{i}(L)| \leq l_0$.

\item $z$ is contained in $\mathbb{C}_{f^{n}(L)}$, $\widehat{(z,f^{n}(J)}\geq \epsilon $ and $dist(z,f^{n}(J)) \leq K(\epsilon,\delta,l_0)|f^{n}(J)| $.

\end{itemize}
  Then the backward orbit of $z$ along  the orbit of $J$ is well defined and

\begin{equation}
\frac{dist(z_{-n},J)}{|J|} \leq C(\delta,\epsilon)\frac{dist(z,f^{n}(J))}{|f^{n}(J)|}
\end{equation}

Moreover, fixing $\epsilon$ and  $\delta$,  $K$ tends to infinity when $l_0$ goes to zero. 
 \end{lem}
\begin{proof} If $z$ is close to $f^{n}(J)$, use lemma \ref{nodist}.i; otherwise, capture $z$ in a Poincare neighborhood whose diameter is commensurable with $dist(z,f^{n}(J))$ and use lemma \ref{pulb}. For details, see \cite{LYA97}. \end{proof}

The proof of the following lemma is easy. Assume that $q$ is a local maximum (other case is similar).

\begin{lem}\label{soum}There exists $\epsilon$ with the following properties. Let $z \in \mathbb{C}$ be a point near to $q$, $q \in \mathbb{I}_p$.

\begin{itemize}

\item If  $\widehat{(f(z),f(q), f(q) + \mathbb{R^{+}})} <  \epsilon$  then $min\{\widehat{(z,q, q + \mathbb{R^{+}})}, \widehat{(z,q, q + \mathbb{R^{-}})}   \} > \epsilon$.

\item If  
\begin{itemize}
\item $\widehat{(f(z),f(q), f(q) + \mathbb{R^{-}})} < \epsilon$. 
\item $min\{\widehat{(z,q, q + \mathbb{R^{+}})}, \widehat{(z,q, q + \mathbb{R^{-}})}\} <  \epsilon$.
\end{itemize}
Then either $f^{+}_q(f(z))=z$ or $f^{-}_q(f(z))=z$. 
\end{itemize} 
\end{lem}

\begin{figure}
\centering
\psfrag{b}{$f$}
\psfrag{ds}{$D_{\alpha}(S^j_{q})$}
\psfrag{dt}{$D(T^j_{q})$}
\psfrag{dvp}{$D_{\tilde{\theta}_q}(V_q)$}
\psfrag{dvq}{$D_{\tilde{\theta}_r}(V_r)$}
\psfrag{du}{$D_{\theta_r}(U_r)$}
\psfrag{pN}[][][0.5]{ $Q_{-(i + N_j)}^{k}$}
\psfrag{pNm}[][][0.5]{ $Q_{-(i + N_j-1)}^{k}$}
\psfrag{pn}[][][0.5]{ $Q_{-(i + n_j^{q})}^{k}$}
\psfrag{pnm}[][][0.5]{ $Q_{-(i + n_j^{q}-1)}^{k}$}
\psfrag{pi}[][][0.5]{$Q_{-i}^{k}$}
\psfrag{d}{$f$}
\psfrag{P}{$Q_{0}^{j}$}
\psfrag{Q}{$R_{0}^{j}$}
\psfrag{a}{$f^{n_{j}^{q}-1}$}
\psfrag{c}{$f^{n_{j}^{r}-1}$}
\includegraphics[width=0.80\textwidth]{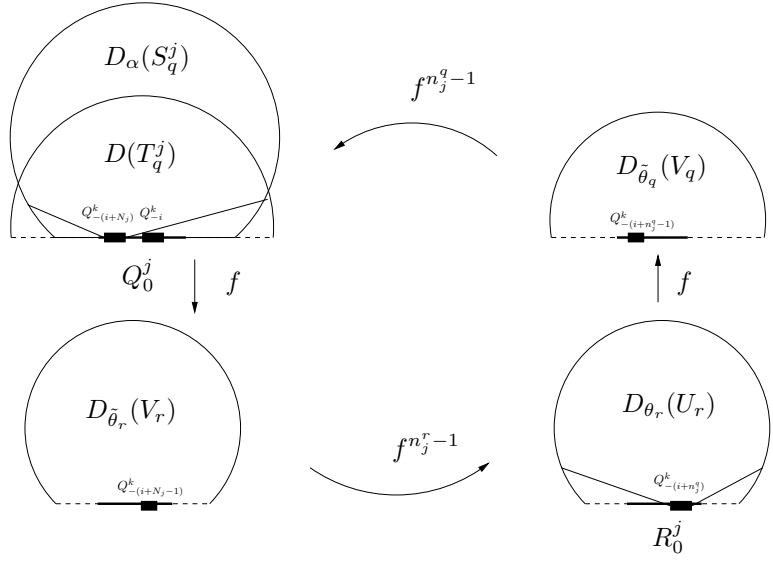}
\caption{  lemma \ref{ere}. In the picture, two critical points, $q$ and $r$, are involved in the renormalization.}
\end{figure}

\section{Complex Bounds}

Let $z \in \mathbb{C}$. We say that $z_{-\ell}$ is well defined if the backward orbit of $z$ along the orbit of $Q_{-\ell}^k$ by $f^\ell$ is well defined.

Recall that the pullback of $T^{k}_{p}$ along the $k$-cycle of renormalization exists. Hence let $U_{p}=T^{k}_{p}, U_{r}, V_{r}, {\tilde{U}}_{p}=S^{k}_{p}$, $r \in \mathbb{I}_p$, be as in definition \ref{pak}. The most important step to prove the complex bounds using the Lyubich-Yampolsky approach is the following lemma:

\begin{mle}[At most linear growth of distances] \label{lg} There exists $C$ with the following property. Provided that $k \geq k_{0}(a)$, if $z \in D(T^{k - a}_q)\cap \mathbb{C}_{U_{q}}$, then $z_{-(n^{q}_{k}-1)}$ is well defined and 

\begin{equation}
\frac{dist(z_{-(n_{k}^{q}-1)},Q_{-(n_{k}^{q} -1)})}{|Q_{-(n_{k}^{q} -1)}|}\leq C\frac{dist(z,Q_{0}^{k})}{|Q_{0}^{k}|} 
\end{equation}
\end{mle}

We say that $\ell$ is a \textbf{return} to $R_{0}^{j}$ if $Q_{-\ell}^{k} \subset R_{0}^{j}$. The time $\ell$ is a \textbf{return to level j} if $\ell$ is a return to $R_{0}^{j}$ for some $r \in \mathbb{I}_{p}$.

In this section, we assume that the levels $j$ and $k$ are deep enough. We say that $z_{-i}$ \textbf{$\epsilon$-jumps} if $\widehat{(z_{-i},Q_{-i}^{k})} > \epsilon$. Fix $a$ and assume that $z \in D(T^{k - a}_q)\setminus \mathbb{R}$ ( for $z \in U_q$, the proof of main lemma follows of the results in section \ref{sect}).

\begin{rem}\label{zeroc} Let $U_q=T_{q}^{j}, U_r, V_r, \tilde{V}_r,  \tilde{U}_q=S_{q}^{j}$ and $\delta > 0$  be as in remark \ref{zeror}. Because $\sum_{i < n^r_j} |f^{i}(\tilde{V}_r)|$ goes exponentially fast to zero, fixing $\theta$, the complex pullback of $D_{\theta}(U_{r})$ until $D_{\tilde{\theta}}(V_{r})$ by $f^{n^r_j-1}$ can be extended to a larger neighborhood $D($$\delta$-$U_r)$,for $j \geq j(\theta)$, by lemma \ref{pulb}. Hence we can use lemma \ref{nodist} with $M=U_r$ and  appropriate inverse branchs $h_{t}$ of $f^{t}$, $t \leq n^{r}_j -1$, such that $h_t(z)=z_{-t}$, where  $z$ is inside $D_\theta($$\delta$-$U_r)$ and $z_{-t}$ is the backward orbit of $z$ along the orbit of $V_r$ by $f^{n^{r}_j -1}$. 
\end{rem}

\begin{lem}[Tour along the j-cycle] \label{ere} There exist $C$ and $\epsilon_0$ with the following property. Assume that  $\ell_0$ is a return to $Q_{0}^{j}$, $j \geq k - a$; $z_{-\ell_0}$ is well defined and $z_{-\ell_0} \in D(T^{j}_{q})$. Let  $\ell_1=min\{n^q_k-1, \ell_0 + N_j\}$. Then $z_{-\ell_1}$ is well defined and
\begin{itemize}
\item[\textbf{a.}] If $\ell$ is a return to $R^{j}_{0}$, for $\ell_0 \leq  \ell \leq \ell_1$, then
$dist(z_{-\ell},Q_{-\ell}^{k}) \leq C |R_{0}^{j}|$.

\item[\textbf{b.}] Let $\epsilon \leq \epsilon_0$. If $\ell_0 + N_j \leq n^q_k-1$ then one of the following statements holds:
\begin{itemize}
\item $z_{-(\ell_0+N_{j})}$ is contained in $D(T^{j}_{q})$.
\item For some return $\ell$ to level j, $\ell_0 \leq \ell \leq \ell_0 + N_{j}$, $z_{-\ell}$ $\epsilon$-jumps.
\end{itemize}

\item[\textbf{c.}] Let $\epsilon \leq \epsilon_0$. Suppose that $z_{-\ell}$ does not $\epsilon$-jump in returns to level $j$, $\ell_0 \leq \ell \leq \ell_1$. Then  $z_{-\ell}$ does not $\alpha$-jump, for each $\ell$ with $\ell_0 \leq \ell < \ell_1$, where $\alpha=O(\epsilon)$.
\end{itemize}
\end{lem}

\begin{proof}(compare with the proof of lemma 5.3 in  \cite{LYA97}) The  pullback of $T^{j}_{q}$ along the $j$-cycle exists. Hence let $U_{q}=T^{j}_{q}, U_{r}, V_{r}, {\tilde{U}}_{q}=S^{j}_{q}$, $r \in \mathbb{I}_p$, be as in definition \ref{pak}. The complex pullback of $D_{\theta_q}(T^{j}_q)$, $\theta_q = \pi/2$,  until $V_{q}$ by $f^{n^q_j -1}$ is inside $D_{\tilde{\theta}_q}(V_{q})$, with $\tilde{\theta}_q \sim \pi/2$, by lemma \ref{pulb}. By remark \ref{zeror} and lemma \ref{cut}, $f(c)$ cuts $V_{q}$ in commensurable parts, where $c$ is the predecessor of $q$ at level $j$. Hence , by lemma \ref{cpulq}, the complex pullback of $D(V_{q})$ by $f$ along the orbit of $c$ is contained in $D_{\theta_{c}}(U_{c})$, where $\theta_c$ does not depend on $j$. By induction, assume that the complex pullback of $D(T^{j}_{q})$ along the $j$-cycle until $U_{r}$ is inside $D_{\theta_{r}}(U_{r})$. Again the complex pullback of $D_{\theta_{r}}(U_{r})$ until $V_{r}$ by $f^{n^r_j - 1}$ is inside $D_{\tilde{\theta}_{r}}(V_{r})$, $\tilde{\theta}_{r} \sim \theta_{r}$, by lemma \ref{pulb}.  The complex pullback of $D_{\theta_r}(V_{r})$ by $f$ along the orbit of $d$ is contained in $D_{\theta_{d}}(U_{d})$, where $d$ is the predecessor of $r$ at level $j$, by remark \ref{zeror}, lemma \ref{cut} and lemma \ref{cpulq}. Hence the complex pullback of $D(T^{j}_{q})$ along all $j$-cycle is inside $D_{\alpha}(S_{q}^{k})$, for some $\alpha$.  Observe that if $Q_{-\ell}^{k} \subset R_{0}^{j}$, for some $r \in \mathbb{I}_p$, then  $z_{-\ell}$ and $Q_{-\ell}^{k}$ are inside $D_{\theta_{r}}(U_{r})$, which  has diameter commensurable with $|R_{0}^{j}|$ (because  $U_{r} \subset R_{0}^{j-1}$). This proves \textbf{a}. 
Assume $\ell_0+N_{j} \leq n^q_k -1$ and suppose that $z_{-(\ell_0+N_{j})}$ is not inside  $D(T^{j}_{q})$. Every point $z \in D_{\alpha}(S^{j}_{q}) \setminus D(T^{j}_{q})$ satisfies $\widehat{(z,J)} \geq \epsilon$, for all interval $J$ inside $S^{j}_{q}$,  for a definitive  $\epsilon$, since $S^{j}_{q}$ is $\delta$-deeply contained in $T^{j}_{q}$. This proves \textbf{b}.
To prove \textbf{c}, suppose that $z_{-\ell}$ does not $\epsilon$-jump in returns to level $j$, $\ell_0 \leq \ell \leq \ell_1$. If $Q^{k}_{-\ell} \subset R^{j}_{0}$ then $z_{-\ell} \in D_{\theta_r}(U_{r})$. By lemma \ref{nodist}.ii and remark \ref{zeroc}, $z_{-i}$ does not $\alpha$-jump for each $\ell \leq i  \leq min\{\ell + n_j^{r}-1,\ell_1\}$, with $\alpha=O(\epsilon)$. 
\end{proof}

The concept of 'good return' was introduced in \cite{LYA97}:

\begin{lem}[Good returns] \label{gr} Assume that $\ell_0$ is a return to $Q_{0}^{j}$, $z_{-\ell_0}$ is well defined and  $z_{-\ell_0} \in D(T^{j}_{q})$. If there exists $\ell$ satisfying 

\begin{itemize}
\item $z_{-\ell}$ is well defined, 
\item $\ell$ is a return to level $j$, $j \geq k - a$, with $\ell_0 \leq \ell \leq N_{j +1}$,
\item $z_\ell$ $\epsilon$-jumps.
\end{itemize}
then provided $k \geq k_0(a)$, $z_{-(n_k^{q}-1)}$ is well defined and 

\begin{equation}
\frac{dist(z_{-(n_k^{q}-1)},Q_{-(n_k^{q}-1)}^{k})}{|Q_{-(n_k^{q}-1)}^{k}|}\leq C(\epsilon)\frac{|Q_{0}^{j}|}{|Q_{0}^{k}|} 
\end{equation}
Furthermore, if $\ell$ is the smallest time satisfying the above properties (in this case we say that $\ell$ is a \textbf{good return} to level $j$)  then $z_{-\ell} \in D_{\theta(\epsilon,k-j)}(Q_{-\ell}^{k})$
\end{lem}

\begin{proof}
 Let $\ell$ be the smallest time satisfying the above properties. We have $Q_{-\ell}^{k} \subset R_{0}^{j}$, for some $r \in \mathbb{I}_{p}$. Then  $dist(z_{-\ell},Q_{-\ell}^{k}) \leq C |R_{0}^{j}|$, by lemma \ref{ere} (use \ref{ere}.b successively, and finally \ref{ere}.a). Since $N_{j+1}/N_{j}$ is bounded,  by lemma \ref{icom} there exists $C$ such that $ |R_{0}^{j}|/|Q_{-\ell}^{k}| \leq C |Q_{0}^{j}|/|Q_{0}^{k}|$. Hence
\begin{equation}
\frac{dist(z_{-\ell},Q_{-\ell}^{k})}{|Q_{-\ell}^{k}|}\leq C\frac{|Q_{0}^{j}|}{|Q_{0}^{k}|}\end{equation} 
Since  $z_{-\ell}$ $\epsilon$-jumps and $|Q_{0}^{j}|/|Q_{0}^{k}| \leq C(k-j)$, by real bounds, $z_{-\ell}$ is contained in  $D_{\theta(k-j,\epsilon)}(Q_{-\ell}^{k})$. By remark \ref{zeror},
there exists an interval of monotonicity $\tilde{V}_q$ for $f^{n_k^{q}-1 - \ell}$ which contains $Q_{-(n_k^{q}-1)}^{k}$ and such that $f^{n_k^{q}-1 - \ell}(L)$ contains $\delta$-$Q_{-\ell}^{k}$, for some $\delta$. Moreover $|Q_{0}^{j}|/|Q_{0}^{k}| \leq K$, where $K$ depends only on $a$. By lemma \ref{liu}

\begin{equation}
\frac{dist(z_{-(n_k^{q}-1)},Q_{-(n_k^{q}-1)}^{k})}{|Q_{-(n_k^{q}-1)}^{k}|} \leq C(\epsilon)
\frac{dist(z_{-\ell},Q_{-\ell}^{k})}{|Q_{-\ell}^{k}|} 
\end{equation}
\end{proof}

\begin{lem}[Inductive step] \label{ind} Let $\alpha$ be as in lemma \ref{ere}.c, $\epsilon$ small enough and $j \geq k-a$ . Provided that $k \geq k_0(a)$, the following holds. Suppose that   $z_{-N_{j}} \in D(T^{j}_{q})$ and $z_{-\ell}$ does not $\alpha$-jump for $\ell < N_{j}$. If $N_{j+1} \leq n^{q}_k -1$, then one of the following statements holds:

\begin{enumerate}
\item There exist $\ell$ with $N_{j} \leq \ell \leq N_{j+1}$ such that 
\begin{itemize}
\item $z_{-\ell}$ is well defined;
\item $\ell$ is a return to level $j$.
\item $z_{-\ell}$ $\epsilon$-jumps.
\end{itemize}
\item $z_{-N_{j+1}}$ is well defined and $z_{-N_{j+1}} \in D(T^{j+1}_{q})$.
\end{enumerate}
\end{lem}

\begin{proof}(compare with the proof of corollary 5.5 in  \cite{LYA97}) Let $L$ be the maximal interval of monotonicity of $f^{N_{j+1}}$ which contains $Q_{-N_{j+1}}^{k}$. Since $Q_{-N_{j+1}}^{k}$ is contained in $Q_{0}^{j+1}$, the interval  $L$ is contained in $L^{j+1}_{q}$. In particular, $\sum_{i \leq N_{j+1}}|f^{i}(L)|$ is very small. Assume that (i) does not hold. By lemma \ref{ere}.c,  $z_{-\ell}$ never  $O(\epsilon)$-jumps for $\ell \leq N_{j+1}$. Hence, by lemma \ref{mi},  $z_{-N_{j+1}} \in D_{(\pi/2 + O(\epsilon))\beta}(L)$, with $\beta \sim 1$. Since  $L^{j+1}_{q}$ is $\delta$-deeply contained in  $T^{j+1}_{q}$,  then $D_{(\pi/2 - O(\epsilon))\beta}(L) \subset D(T^{j+1}_{q})$, if  $\epsilon$ is small enough. \end{proof}

For a deep level $\ell$, let $\mathbb{I}_p=\{c^{\ell}_i\}_{i \leq \#\mathbb{I}_p}$, where  $c^\ell_1=q$ and $c^\ell_i$ is the predecessor of $c^\ell_{i-1}$ at level $\ell$.

\begin{proof}[Proof of \textbf{Main Lemma \ref{lg}}](compare \cite{LYA97}, subsection 5.2)
Let $b$, with  $k - a \leq b \leq k$, be such that  $z \in D(T^{b}_{q})\setminus D(T^{b+1}_{q})$. In other words, $z$ is in the 'scale' of level  $b$. If $z \in D(U_q)$, we are done by lemma \ref{nodist}.i and remark \ref{zeroc}. Otherwise $dist(z,Q_{0}^{k})/|Q_{0}^{k}| \sim |Q_{0}^{b}|/|Q_{0}^{k}|$. Note that $n_k^{q} -1 = CN_b + \sum_{t < \tilde{t}} n_b^{c^b_t} + (n_b^{r} -1)$, where $r=c^b_{\tilde{t}}$ and  $\tilde{t} \leq \# \mathbb{I}_p$. If $C$ is zero, using arguments as in the proof of lemma \ref{ere}.a, $z_{-(n_k^{q} -n_b^r) }$ is well defined and it is contained in $D_{\theta_r}(U_r)$, where $U_r$ and $\theta_r$ are as in the proof of \ref{ere}.a (replacing $j$ by $b$). In particular $dist(z_{-(n_k^{q} -n_b^r)},Q^k_{-(n_k^{q} -n_b^r)}) \leq C |R^b_{0}|$. Since $n_k^{q} -n_b^r \leq N_{b+1}$, by lemma \ref{icom} $|R^b_{0}|/|Q^k_{-(n_k^{q} -n_b^r)}| \leq C |Q^b_{0}|/|Q^k_{0}|$. By lemma \ref{nodist}.i(recall remark \ref{zeroc}), we obtain the lemma. Otherwise, suppose that $z_\ell$ $\epsilon$-jumps for some $\ell \leq N_b$, where  $\ell$ is a return to level $b$. Then we are done, by lemma \ref{gr}. If $\epsilon$-jumps in a return to level $b$ does not  happen then $z_{-N_b} \in D(T^{b}_{q})$, by lemma \ref{ere}.b. Let $i$ be the deepest level such that $n_k^{q} -1 = CN_i + \sum_{t < \widehat{t}} n_i^{c^i_t} + (n_i^{r} -1)$,  with  $C$ distinct of zero, $r=c^i_{\widehat{t}}$ and  $\widehat{t} \leq \# \mathbb{I}_p$. Note that $C$ is bounded. Suppose, by induction, that $z_{-N_j} \in D(T^{j}_{q})$, $j < i$, and $z_\ell$ does not $\alpha$-jump, for $\ell \leq  N_j$. By lemma \ref{ind}, either there is  a return to level $j$ between $N_j$ and $N_{j+1}$, which $\epsilon$-jumps, or $z_{-N_{j+1}} \in D(T^{j+1}_{q})$ and moreover $z_\ell$ does not $O(\epsilon)$-jump for $\ell \leq N_{j+1}$, by lemma \ref{ere}.c. If there are not good returns between the times $0$ and $N_i$ then $z_{-N_i} \in D(T^i_{q})$. Then either there are good returns to level $i$ between  $N_i$ and $CN_i$, and we can use lemma \ref{gr}, or $z_{-CN_i} \in D(T^i_{q})$. Using arguments as in the proof of \ref{ere}.a, $z_{-(n_k^{q}  -n_i^r)}$ is well defined and it is inside $D_{\theta_r}(U_r)$, where $U_r$ and $\theta_r$ are as in the proof of lemma \ref{ere}.a (replacing $j$ by $i$). In particular $dist(z_{-(n_k^{q}  -n_i^r)},Q^k_{-(n_k^{q} -n_i^r)}) \leq C |R^i_{0}|$. Since $n_k^{q} -n_i^r \leq N_{i+1}$, one gets $|R^i_{0}|/|Q^k_{-(n_k^{q} -n_i^r)}| \leq C |Q^i_{0}|/|Q^k_{0}|$, by lemma \ref{icom}. By lemma \ref{nodist}.i(recall remark \ref{zeroc}), we obtain the lemma.\end{proof}

\begin{defi} We say that $z \in \mathbb{C}$ has \textbf{complex pullbacks along the $k$-cycle} if given, for each $q \in \mathbb{I}_p$, arbitrary inverse branches $g_q$ defined near to $f(q)$, there exists a sequence $z_{-\ell}$, $\ell \leq N_{k}$ which satisfies
\begin{itemize}
\item If $f^{N_k -\ell}(P^{k}_0) \subset Q^{k}_{0}$, $t \leq n^q_k - 1$, then $z_{-(\ell+t)}$, $t < n^q_k - 1$, is the backward orbit  of $z_{-\ell}$ along the orbit of $Q^{k}_{-(n^q_k - 1)}$ by $f^{n^q_k - 1}$.

\item If $f^{N_k -\ell}(P^{k}_0) \subset Q^{k}_{-(n^q_k -1)}$, then  $z_{-\ell}$ is contained in the domain of $g_r$ and  $g_r(z_{-\ell})=z_{-(\ell+1)} $. Here $r$ is the predecessor of $q$ at level $k$.
\end{itemize}
\end{defi}

The following lemma say that a point which never $\epsilon$-jumps has a real itinerary:

\begin{lem} \label{neverj} Assume that $z$ has complex pullbacks along the $k$-cycle. Let $z_{-\ell}$ be a complex pullback such that 
\begin{itemize}

\item If $f^{N_k -\ell}(P^{k}_0) \subset Q^{k}_{-j}$, for some $q \in \mathbb{I}_p$, then $\widehat{(z_{-\ell},Q^{k}_{-j})} \leq \epsilon$, for $\epsilon$ small enough.

\item If $f^{N_k -\ell}(P^{k}_0) \subset Q^{k}_{0}$, for some $q \in \mathbb{I}_p$, Then  $dist(z_{-\ell},Q^{k}_{0}) \leq C  |Q^{k}_{0}|$.
\end{itemize}

There exists $k_0(C)$ such that if $k \geq k_0$ then $z_{-N_k}$ has the itinerary of a point in $\partial P_{0}^{k}$ by $f^{N_k}$.
\end{lem}

\begin{proof} Consider $i$ such that $f^{N_k -i}(P^{k}_0)$ is contained in $Q^{k}_{-(n^{q}-1)}$. Since $\widehat{(z_{-i},Q^{k}_{-(n^{q}_k-1)})} \leq \epsilon$, $z_{-i}$ is inside a $\epsilon$-sector of angle $\epsilon$ which has as its vertex a point in $\partial Q^{k}_{-(n^{q}_k-1)}$. We claim that this vertex is the point $b$ such that $\{ b \}= f(\partial R_{0}^{k})$, where $r$ is the predecessor of $q$ at level $k$. Otherwise, $z_{-(i+1)}$  $\epsilon$-jumps, by the first statement in lemma \ref{soum}. Moreover, by the second statement in lemma \ref{soum}, the inverse branch $g$, defined near to $f(r)$, such that $g(z_{-i})= z_{-(i+1)}$ must be real. Let $R_{0}^{k} = [x_1,x_2]$. By lemma \ref{naomuda}, if $z_{-(i+1)}$ is contained in the $\epsilon$-sector whose vertex is $x_a$, then $z_{-(i+1+t)}$, $t < n_k^{r}$, is contained in the $\epsilon$-sector whose vertex $\tilde{x}_a$ satisfies $f^{t}(\tilde{x}_a)=x_a$. Using the previous arguments along all $k$-cycle, the lemma follows.
\end{proof}

\begin{lem} \label{quad} Let $f(z) = z^{2n}$.  Let $P$ be a symmetric interval $M$ such that  $f(P) \subset M$, with $|M| \leq C|f(P)|$. Then there exist $C_{1}(C,n)$, $C_{2}(C,n)$  such that

\begin{equation}
\frac{dist(z,P)}{|P|} \leq C_{1} (\frac{dist(f(z),M)}{|M|})^{\frac{1}{2n}} + C_{2}
\end{equation}
\end{lem}

\begin{defi} \label{plike} The renormalization $R^{k}_p$ has a \textbf{polynomial-like extension} if there exist  simply connected domains $\tilde{D}_p$ and $D_q$, for each $q \in \mathbb{I}_p$, satisfying

\begin{itemize}
 \item The critical point $q$  is the unique critical point of $f$ in $D_q$. Moreover, if $q'_k \neq p$, $f^{n^{q'_k}_k} \colon D_q \rightarrow D_{q'_k}$ is a ramified covering.
 
 \item If $q$ is the predecessor of $p$ at level $k$, then $f^{n^{p}_k} \colon D_q \rightarrow \tilde{D}_p$ is a ramified covering. Moreover the closure of $D_p$ is contained in the interior of $\tilde{D}_p$. 
\end{itemize}

In particular $f^{N_k}\colon D_q \rightarrow \tilde{D}_p$ is a polynomial-like map.
\end{defi}

The $\alpha$-sector supported on a interval $J$ is the set of complex number which satisfies  $\widehat{(z,P_{0}^{k})} \geq \alpha$.

\begin{proof}[Proof of Theorem \ref{cb}](compare \cite{LYA97}, subsection 5.3) By the Main Lemma and  lemma \ref{quad}, if  $dist(z,P_{0}^{k})/|P_{0}^{k}| \leq C$ then $z \in \mathbb{C} \setminus \mathbb{R}$ has complex pullbacks along the $k$-cycle. Moreover, if $\ell_r$ is such that $f^{N_{k}-\ell_{r}}(P_{0}^{k})$ is contained in $Q_{0}^{k}$, then $dist(z_{-\ell_r}, Q^{k}_{0}) \leq \tilde{C}|Q^{k}_{0}|$. Hence $z_{-\ell_r} \in D(T^{a_k}_r)$, where $k - a_k$ is bounded. Furthermore,  for  points satisfying $dist(z,P_{0}^{k})/|P_{0}^{k}|= C$, since $C$ is large enough, one has

\begin{equation} \label{contr}
\frac{dist(z_{-N_{k}},P_{0}^{k})}{|P_{0}^{k}|} \leq \frac{1}{10} \frac{dist(z,P_{0}^{k})}{|P_{0}^{k}|}
\end{equation}
Hence, let  $D$ be a disc large enough such that (\ref{contr}) holds in $\partial D$. Then the pullback of $\partial D \cap (\mathbb{C} \setminus \mathbb{R})$ along the $k$-cycle  is contained in $D \setminus \mathbb{R}$. Using arguments as in the proof of lemma \ref{ere}.a, the complex pullbacks of $T^{k}_p$ are inside $D_{\theta}(S^{k}_p)$. Hence, taking a  larger disc, if necessary, the pullback of $V = D\cap \mathbb{C}_{T^{k}_p}$ is compactly contained in $V \cap \mathbb{C}_{L^{k}_p}$ (clearly the pullback is contained in $V \cap \mathbb{C}_{S^{k}_p}$, but this could not be true for the closure of the pullback).  One gets a polynomial-like  extension. Now, it is necessary to control the modulus of the polynomial-like extension. By lemma 2.4 in  \cite{LYA97}, it is sufficient  to prove that the ring $V \setminus K(R^{k}_p(f))$ has the  modulus bounded below, where $K(R^{k}_p(f))$ is the filled-in Julia set of $R^{k}_p(f)$. We claim that $K(R^{k}_p(f))$ is contained in a sector supported on $L^{k}_p$. Indeed, take $z' \in K(R^{k}_p(f))$. Let $z = f^{N_{k}}(z')$. If $z \in T^{k}_p$, one gets that $z' \in D_{\theta}(S^{k}_p)$. Assume that $z \in \mathbb{C}\setminus \mathbb{R}$. Denote $z_{-\ell} = f^{N_k -\ell}(z')$. If there exist $t_0 < t$, $j$ and $r \in \mathbb{I}_p$ satisfying the following properties
\begin{itemize}
\item $R^k_{-t_0}$ is contained in $R^{j}_0$.
\item $z_{-(\ell_r + t_0)}$ is contained in $D(T^j_r)$.
\item $R^k_{-t}$ is  contained in $Q^j_0$, for some $q \in \mathbb{I}_p$. Moreover $ t \leq N_{j+1}$.
\item $\widehat{(z_{-(\ell_r +t )},R^k_{-t})\geq \epsilon}$.  
\item $t$ is minimal with the above properties.
\end{itemize}
then we capture  $z_{-(\ell_r + t )}$ in  $D_{\gamma}(R_{-t}^{k})$, by lemma \ref{gr}. Using arguments as in the proof of lemma \ref{ere}.a, we obtain that $z' \in D_{\tilde{\gamma}}(P_{0}^{k})$. In the other case, $z_{-\ell}$ never $O(\epsilon)$-jumps for all $\ell$. Then, by lemma \ref{neverj}, $z'$ has the itinerary of a point in $\partial P_{0}^{k}$. Let $L$ be the interval of monotonicity of $f^{N_k}$ which contains this point. Clearly $L \subset L^{k}_{p} \subset_{\delta} T^{k}_{p}$. Since $L^{k}_{p}$ is contained in $P^{k-1}_0$ and $N_k/N_{k-1}$ is bounded, $\sum_{i \leq N_k} |f^i(L)|$ goes exponentially fast to zero. Hence By lemma \ref{mi}, $z' \in D_{(\pi/2 +O(\epsilon))\beta}(L)$, where $\beta \sim 1$, by lemma \ref{mi}. We captured the Julia set in three Poincare neighborhoods: two of these neighborhoods are based on intervals contained in $L^{k}_p$ which contain points in $\partial P_{0}^{k}$. The other one is $D_{max\{\tilde{\gamma},\theta \}}(L_{p}^{k})$. Hence the Julia set is contained in $D_{\theta}(L^{k}_{p})$, for some angle $\theta$.
\begin{figure}
\centering
\psfrag{r}{$u$}
\psfrag{r'}{$\tilde{u}$}
\psfrag{p}{$p$}
\includegraphics[width=0.50\textwidth]{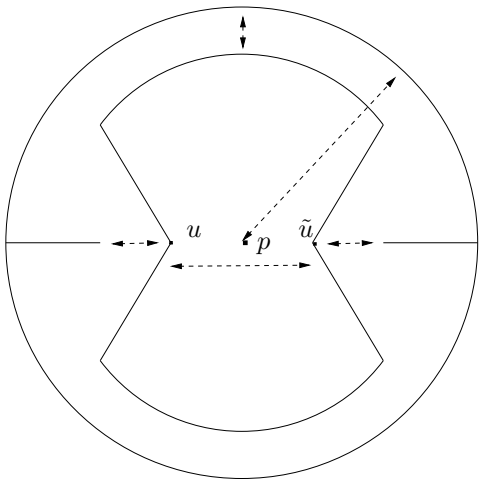}
\caption{ \label{julia} This ring is contained in $V \setminus K(R^{k}_p(f))$. The extern boundary is $\partial V$. All measures indicated in picture are commensurable. Moreover, the angles in $u$ and $\tilde{u}$ are under control.}
\end{figure}
Indeed, the filled-in Julia set is inside a sector supported on $P_{0}^{k}$. We follows the argument in \cite{LYA97}: Consider the interval $J$ of monotonicity of $f^{N_{k}}$ which contains the periodic point $u$ in $\partial P_{0}^{k}$. We have $J \subset_{\delta} f^{N_{k}}(J)$. Let $g = f^{-N_k}$ be the inverse branch of $f^{N_{k}}$ such that $g(D(f^{N_{k}}(J)))$ is deeply inside $D(J)$, for large $k$. The map $g$ has exactly one fixed point $u$ in $D(J)$ and $\lambda_1 \geq (f^{N_k})'(u) \geq \lambda_2 > 1$, where $\lambda_1$, $\lambda_2$ does not depend on $k$. We can linearize $g$ inside $D(f^{N_{k}}(J))$. Do the pullback of $D_{\theta}(L^{k}_{p})$ along the k-cycle a sufficient number of times (which does not depend on $k$) such that the pullback will be contained in $D_{\tilde{\theta}}(H)$, for some interval $H$ such that the  component of $H\setminus P_{0}^{k}$ which contains $u$ is contained in $g^2(B_{\delta|P_0^k|}(u))$, for a small $\delta$. In particular, since $K(R^{k}_p(f)) \subset D_{\tilde{\theta}}(H)$,  the Julia set cuts the fundamental ring $B_{\delta|P_0^k|}(u)\setminus g(B_{\delta|P_0^k|}(u))$ at most a certain angle. Using that the linearization has small distortion on $B_{\delta|P_0^k|}(u)$, if $\delta$ is small, and the invariance of $K(R^{k}_p(f))$ by $R^{k}_p(f)$, we obtain the sector. By invariance of $K(R^{k}_p(f))$, we obtain the sector in other point of $\partial P_{0}^{k}$.
\end{proof}

\begin{rem} If $f$ satisfies the standard conditions, the complex extensions can be taken unbranched. In fact, there exists $\delta$ such that $\delta$-$P_{0}^{k} \cap P(f) = P_{0}^{k} \cap P(f)$. By the control of geometry of the ring $V\setminus K(R^k_p)$ (see the figure \ref{julia}), the ring $(V\setminus K(R^k_p(f))) \cap$ $\mathbb{C}_{\delta\text{-} P_{0}^{k}}$ has the modulus bounded below. Now it is easy to find an unbranched polynomial like extension such that the modulus is also bounded below.
\end{rem}

\section{Acknowledgments.} I wish to thank W. de Melo for introducing me on this subject and for useful conversations about the mathematics and the style of this article. I am also grateful to A. Pinto for comments that suggested lemma  \ref{icom}. This work was supported by CNPq-Brazil .

\end{document}